# A few equalities involving integrals of the logarithm of the Riemann $\varsigma$-function and equivalent to the Riemann hypothesis II

Sergey K. Sekatskii, Stefano Beltraminelli, and Danilo Merlini

**Abstract.** This paper is a continuation of our recent paper with the same title, arXiv:0806.1596v1 [math.NT], where a number of integral equalities involving integrals of the logarithm of the Riemann $\varsigma$-function were introduced and it was shown that some of them are equivalent to the Riemann hypothesis. A few new equalities of this type are established; contrary to the preceding paper the focus now is on integrals involving the argument of the Riemann function (imaginary part of logarithm) rather than the logarithm of its module (real part of logarithm). Preliminary results of the numerical research performed using these equalities to test the Riemann hypothesis are presented. Our integral equalities, together with the equalities given in the previous paper, include all earlier known criteria of this kind, viz. Wang, Volchkov and Balazard-Saias-Yor criteria, which are certain particular cases of the general approach proposed.

## 1. Introduction

In a recent paper [1] we analyzed integrals of the type $\int_0^\infty \frac{\ln|\varsigma(b+it)|}{a^2+t^2}dt$ (*a, b* real positive) and established a number of equalities equivalent to the Riemann Hypothesis (RH; here $\varsigma(z)$ is Riemann zeta-function, see e.g. [2] for definitions and discussion of the general properties of this function). Our main idea was to introduce the function $g(z)=\frac{1}{a^2-(z-b)^2}$ and consider a contour integral $\int_C \ln(\varsigma(z))g(z)dz$ taken around the rectangular contour *C* with the vertices $b-iX$, $b+iX$, $b+X+iX$, $b+X-iX$. On the left border of the contour, that is on the line *z=b+it*, this integral takes the required form $-i\int_{-X}^{X}\frac{\ln(\varsigma(b+it))}{a^2+t^2}dt$ and it is not difficult to demonstrate [1, 3] that the integrals taken along the external lines of the contour tend to zero when $X\to\infty$. Then the analysis of the contour integrals using the residue and (generalized) Littlewood theorems [1, 4] give required results.

In this short note we analyze, along the same lines, the integrals $\int_0^\infty \frac{t\arg(\varsigma(b+it))}{(c^2+t^2)(d^2+t^2)}dt$, $\int_0^\infty \frac{t\arg(\varsigma(b+it))}{(a^2+t^2)^2}dt$, which involve an argument of the Riemann $\varsigma$–function rather than the logarithm of its module, and establish a few other equalities equivalent to the RH.

## 1. Integral equalities involving an argument of the Riemann $\varsigma$ - function

*2.1 Integrals* $\int_0^\infty \frac{t\arg(\varsigma(b+it))}{(c^2+t^2)(d^2+t^2)}dt$

Let us consider the same contour *C*, introduce the function $g(z)=-i\frac{z-b}{(c^2-(z-b)^2)(d^2-(z-b)^2)}$, where *b>-2* is real (with this choice we simply want to avoid unnecessary complications with the trivial Riemann zeroes), $c\ne d$ real positive, and consider the contour integral $\int_C \ln(\varsigma(z))g(z)dz$. Again, the disappearance of the integral taken along the external lines of the contour is certain while on the left border of the contour, that is on the line *z=b+it*, it takes the form $-i\int_{-\infty}^{\infty}\frac{t\ln(\varsigma(b+it))}{(c^2+t^2)(d^2+t^2)}dt$. Function $g(t)=\frac{t}{(c^2+t^2)(d^2+t^2)}$ is an odd function of *t* as well as the function $\arg(\varsigma(b+it))$ while the function $\ln|(\varsigma(b+it)|$ is even. Thus we have that the contour integral now reduces to the form $2\int_0^\infty \frac{t\arg(\varsigma(b+it))}{(c^2+t^2)(d^2+t^2)}dt$. Throughout the paper, similarly to our earlier paper [1], we consider the value of $\arg(\varsigma(b+it))$ obtained by continuous variation along the straight lines joining the points *2, 2+it, b+it* starting from the value 0, cf. [2] P. 212.

In the interior of the contour we have simple poles at $z=b+c$, $z=b+d$ and, if *b<1*, also a single pole of the Riemann function at *z=1*. If *b<1*, in the interior of the contour we also can have a number of zeroes of the Riemann function, and we definitely have an infinite number of them if *b<1/2*. The contribution of the poles of *g(z)* to the contour integral value is, by a residue theorem, equal to

$$2\pi i(i\frac{1}{2(c^2-d^2)}\ln|\varsigma(b+d)|+i\frac{1}{2(d^2-c^2)}\ln|\varsigma(b+c)|) = \frac{\pi}{(d^2-c^2)}\ln|\frac{\varsigma(b+d)}{\varsigma(b+c)}|.$$

The corresponding contribution of the pole of the Riemann function is, according to the generalized Littlewood theorem [1], for *b<1* and $b+c\neq 1$, $b+d\neq 1$, equal to

$$2\pi i\int_b^1 -i\frac{z-b}{(c^2-(z-b)^2)(d^2-(z-b)^2)}dz = 2\pi\int_0^{1-b}\frac{p}{(p^2-c^2)(p^2-d^2)}dp =$$

$$\frac{\pi}{c^2-d^2}\ln|\frac{p^2-c^2}{p^2-d^2}|\Big|_0^{1-b} = \frac{\pi}{d^2-c^2}\ln|\frac{(d^2-(1-b)^2)c^2}{(c^2-(1-b)^2)d^2}|$$

Analogously, the order *n* zero of the Riemann function $\rho=\sigma_k+it_k$, for $b<\sigma_k$, contributes

$$\frac{\pi n}{c^2-d^2}\ln|\frac{d^2-p^2}{c^2-p^2}|\Big|_{it_k}^{\sigma_k-b+it_k} = \frac{\pi n}{c^2-d^2}(\ln|\frac{d^2-(\sigma_k-b)^2+t_k^2-2it_k(\sigma_k-b)}{c^2-(\sigma_k-b)^2+t_k^2-2it_k(\sigma_k-b)}|-\ln|\frac{d^2+t_k^2}{c^2+t_k^2}|)$$

to the contour integral value.

Collecting everything together, we have the following equality: for *-2<b<1*

$$\int_0^\infty \frac{t\arg(\varsigma(b+it))}{(c^2+t^2)(d^2+t^2)}dt = \frac{\pi}{2(d^2-c^2)}\ln|\frac{\varsigma(b+d)}{\varsigma(b+c)}| + \frac{\pi}{2(d^2-c^2)}\ln|\frac{(d^2-(1-b)^2)c^2}{(c^2-(1-b)^2)d^2}| -$$

$$\frac{\pi}{2(d^2-c^2)}\sum_{\rho,\sigma_k>b,t_k>0}(\ln|\frac{(d^2-(\sigma_k-b)^2+t_k^2)^2+4t_k^2(\sigma_k-b)^2}{(c^2-(\sigma_k-b)^2+t_k^2)^2+4t_k^2(\sigma_k-b)^2}|-2\ln|\frac{d^2+t_k^2}{c^2+t_k^2}|) \quad (1).$$

The sum here is taken over all zeroes $\rho=\sigma_k+it_k$ with $b<\sigma_k$ taken into account their multiplicities; when obtaining (1) we paired the contributions of the complex conjugate zeroes $\rho=\sigma_k\pm it_k$ hence $t_k>0$. The convergence of this sum is guaranteed by the well known properties of the Riemann function zeroes distribution [2]. Of course, for $b>1$ we have simply
$$\int_0^\infty \frac{t\arg(\varsigma(b+it))}{(c^2+t^2)(d^2+t^2)}dt = \frac{\pi}{2(d^2-c^2)}\ln\left|\frac{\varsigma(b+d)}{\varsigma(b+c)}\right|.$$

Now let us take $b \geq 1/2$ and rewrite the contribution of the Riemann function zero $I_\rho = -2\pi \int_{it_k}^{\sigma_k-b+it_k} \frac{pdp}{(c^2-p^2)(d^2-p^2)}$ in somewhat different form. We can use the variable change $q = p+it_k$ to get
$$I_\rho/2\pi = -\int_0^{\sigma_k-b} \frac{q+it_k}{(c^2-q^2+t_k^2-2iqt_k)(d^2-q^2+t_k^2-2iqt_k)}dq$$
from which it immediately follows that the real part of this integral is an integral of the expression
$$-\frac{q(c^2-q^2+t_k^2)(d^2-q^2+t_k^2)-4q^3t_k^2-2qt_k^2(d^2-q^2+t_k^2+c^2-q^2+t_k^2)}{((c^2-q^2+t_k^2)^2+4q^2t_k^2)((d^2-q^2+t_k^2)^2+4q^2t_k^2)}=$$
$$-q\frac{-3t_k^4-2t_k^2(q^2+c^2+d^2)-q^2(c^2+d^2)+c^2d^2+q^4}{((c^2-q^2+t_k^2)^2+4q^2t_k^2)((d^2-q^2+t_k^2)^2+4q^2t_k^2)}.$$

By our choice of $b$ we have $0 \leq q < 1/2$ (all Riemann zeroes lye in the critical strip), and this is known that Riemann zeroes with $\sigma_k > 1/2$, if they exist, have very large values of the imaginary part. Hence, provided $c$ and $d$ are not very large, the sign of this expression is certainly positive: it suffices that $c,d \leq \sqrt{3}t_{k,\min}$. It is announced by X. Gourdon and P. Demichel that the first $10^{13}$ Riemann zeroes have $\sigma_k = 1/2$ (see http://numbers.computation.free.fr/Constants/Miscellaneous/zetazeros1e13-1e24.pdf), but apparently a corresponding paper was never properly published. Based on Ref. [5], there are no such zeroes for $t<545,439,823.215$ hence we can safely take $c,d \leq 9.4 \cdot 10^8$ and formulate the following

THEOREM 1. An equality
$$\int_0^\infty \frac{t\arg(\varsigma(b+it))}{(c^2+t^2)(d^2+t^2)}dt = \frac{\pi}{2(d^2-c^2)}\ln\left|\frac{\varsigma(b+d)}{\varsigma(b+c)}\right| + \frac{\pi}{2(d^2-c^2)}\ln\left|\frac{(d^2-(1-b)^2)c^2}{(c^2-(1-b)^2)d^2}\right| \quad (2),$$
where $b$, $c$, $d$ are real positive numbers such that $c,d \leq 9.4 \cdot 10^8$, $c \neq d$, $b+c \neq 1$, $b+d \neq 1$ and $1 > b \geq 1/2$, holds true for some $b$ if and only if there are no Riemann function zeroes with $\sigma > b$. For $b=1/2$ this equality is equivalent to the Riemann hypothesis.

As an example one can take $b=1/2$, $c=3/2$, $d=7/2$ and then in view of $\varsigma(2) = \pi^2/6$, $\varsigma(4) = \pi^4/90$ obtain the following rather elegant criterion equivalent to RH:

$$\int_0^\infty \frac{t\arg(\varsigma(1/2+it))}{(9/4+t^2)(49/4+t^2)} dt = \frac{\pi}{20}\ln(\frac{\varsigma(4)}{\varsigma(2)}) + \frac{\pi}{20}\ln(54/49) = \frac{\pi}{20}\ln(18\pi^2/245) \quad (3).$$

This equality has been tested numerically, see below.

The case $b+c=1$ or $b+d=1$ actually does not create any problem because the corresponding limit is quite transparent. Let $b+d=1$. After the substitution $b+d = 1+\delta$ and, correspondingly, $d^2 - (1-b)^2 = \delta(2-2b+\delta)$ into

$$\frac{\pi}{2(d^2-c^2)}\ln\left|\frac{\varsigma(b+d)}{\varsigma(b+c)}\right| + \frac{\pi}{2(d^2-c^2)}\ln\left|\frac{(d^2-(1-b)^2)c^2}{(c^2-(1-b)^2)d^2}\right|$$ and rearrangement, it takes the form

$$-\frac{\pi}{2(d^2-c^2)}\ln(|\varsigma(b+c)\cdot(c^2-(1-b)^2)d^2/c^2|) + \frac{\pi}{2(d^2-c^2)}\ln(|\varsigma(1+\delta)\cdot\delta(2-2b+\delta)|)$$

In the first term one can simply put $d=1-b$. For the second term, using the Laurent expansion of the Riemann function in the vicinity of $z=1$: $\varsigma(1+\delta) = \frac{1}{\delta} + \gamma + ...$ ($\gamma$ is the Euler constant), we have that for $\delta \to 0$ $\ln(|\varsigma(1+\delta)\cdot(2\delta(1-b)+\delta^2)|)$ tends to $\ln(2(1-b))$ hence

$$\int_0^\infty \frac{t\arg(\varsigma(b+it))}{(c^2+t^2)((1-b)^2+t^2)} dt = \frac{\pi}{2(c^2-(1-b)^2)}\ln\left|\frac{\varsigma(b+c)(c^2-(1-b)^2)(1-b)}{2c^2}\right| +$$

$$\frac{\pi}{2(c^2-(1-b)^2)}\sum_{\rho,\sigma_k>b,t_k>0}(\ln\frac{((1-b)^2-(\sigma_k-b)^2+t_k^2)^2 + 4t_k^2(\sigma_k-b)^2}{(c^2-(\sigma_k-b)^2+t_k^2)^2 + 4t_k^2(\sigma_k-b)^2} - 2\ln\frac{(1-b)^2+t_k^2}{c^2+t_k^2}) \quad (4).$$

The sum is again over all zeroes $\rho = \sigma_k + it_k$ with $b < \sigma_k$ taken into account their multiplicities.

Analogously as above, we have

THEOREM 1a. An equality

$$\int_0^\infty \frac{t\arg(\varsigma(b+it))}{(c^2+t^2)((1-b)^2+t^2)} dt = \frac{\pi}{2(c^2-(1-b)^2)}\ln\left|\frac{\varsigma(b+c)(c^2-(1-b)^2)(1-b)}{2c^2}\right| \quad (5)$$

where $b$, $c$ are real positive numbers such that $c \leq 9.4 \cdot 10^8$, $b+c \neq 1$ and $1 > b \geq 1/2$, holds true for some $b$ if and only if there are no Riemann function zeroes with $\sigma > b$. For $b=1/2$ this equality is equivalent to the Riemann hypothesis.

As an illustration, we can take, for example, $b=1/2$, $c=3/2$ and then in view of $\varsigma(2) = \pi^2/6$ obtain another rather elegant equality equivalent to RH:

$$\int_0^\infty \frac{t\arg(\varsigma(1/2+it))}{(9/4+t^2)(1/4+t^2)} dt = \frac{\pi}{4}\ln(\frac{\pi^2}{27}) \quad (6)$$

This equality also has been tested numerically.

## 2.2 Integrals $\int\limits_0^\infty \frac{t\arg(\varsigma(b+it))}{(a^2+t^2)^2}dt$

The case of the integrals $\int\limits_0^\infty \frac{t\arg(\varsigma(b+it))}{(a^2+t^2)^2}dt$ is considered similarly introducing the function $g(z) = -i\frac{z-b}{(a^2-(z-b)^2)^2}$ and analyzing the contour integral $\int_C \ln(\varsigma(z))g(z)dz$ taken round the same contour $C$. Now in the interior of the contour we have the pole of the second order at $z=a+b$ which, according to the residue theorem, contributes $2\pi i \frac{d}{dz}\left(-i\frac{z-b}{(a+z-b)^2}\ln(\varsigma(z))\right)\Big|_{z=b+a} = \frac{\pi}{2a}\frac{\varsigma'(a+b)}{\varsigma(a+b)}$ to the contour integral value. For $b<1$, $a+b \neq 1$ we have also the contribution of the simple pole of the Riemann function at the point $z=1$ equal to

$$2\pi i \int\limits_b^1 -i\frac{z-b}{(a^2-(z-b)^2)^2}dz = 2\pi\int\limits_0^{1-b}\frac{p}{(a^2-p^2)^2}dp = \pi\frac{1}{a^2-p^2}\Big|_0^{1-b} = \pi\left(\frac{1}{a^2-(1-b)^2}-\frac{1}{a^2}\right).$$

Similarly, the contribution of an order $n$ Riemann zero $\rho = \sigma_k + it_k$ for $b < \sigma_k$ is equal to $I_\rho = -\pi\frac{n}{a^2-p^2}\Big|_{it_k}^{\sigma_k-b+it_k} = -\pi n\left(\frac{1}{a^2-(\sigma_k-b+it_k)^2} - \frac{1}{a^2+t_k^2}\right)$. Pairing the contributions of complex conjugate zeroes we have

$$I_\rho + I_{\rho^*} = 2\pi n\left(\frac{1}{a^2+t_k^2} - \frac{a^2-(\sigma_k-b)^2+t_k^2}{(a^2-(\sigma_k-b)^2+t_k^2)^2+4t_k^2(\sigma_k-b)^2}\right) =$$

$$2n\pi\frac{(\sigma_k-b)^2(3t_k^2-a^2+(\sigma_k-b)^2)}{(a^2+t_k^2)[(a^2-(\sigma_k-b)^2+t_k^2)^2+4t_k^2(\sigma_k-b)^2]}$$

Again, taking $a \leq \sqrt{3}t_{k,\min}$ suffices to ensure the same signs of all Riemann function zero contributions.

Thus, combining everything together, we obtain the following equality: for $-2<b<1$, $a+b \neq 1$

$$\int\limits_0^\infty \frac{t\arg(\varsigma(b+it))}{(a^2+t^2)^2}dt = \frac{\pi}{4a}\frac{\varsigma'(a+b)}{\varsigma(a+b)} + \frac{\pi}{2}\left(\frac{1}{a^2-(1-b)^2}-\frac{1}{a^2}\right) + \quad (7)$$

$$\pi\sum_{\rho,\sigma_k>b,t_k>0}\left(\frac{1}{a^2+t_k^2} - \frac{a^2-(\sigma_k-b)^2+t_k^2}{(a^2-(\sigma_k-b)^2+t_k^2)^2+4t_k^2(\sigma_k-b)^2}\right)$$

and the following

THEOREM 2. An equality

$$\int_0^\infty \frac{t\arg(\varsigma(b+it))}{(a^2+t^2)^2}dt = \frac{\pi}{4a}\frac{\varsigma'(a+b)}{\varsigma(a+b)} + \frac{\pi}{2}(\frac{1}{a^2-(1-b)^2} - \frac{1}{a^2}) \qquad (8)$$

where $a$, $b$ are real positive, $a+b \neq 1$, $a \leq 9.4 \cdot 10^8$ and $1 > b \geq 1/2$, holds true for some $b$ if and only if there are no Riemann function zeroes with $\sigma > b$. For $b=1/2$ this equality is equivalent to the Riemann hypothesis.

Again, the case of $a+b=1$ is a transparent limit of above expressions for $a+b=1+\delta$, $\delta \to 0$. Using the same Laurent expansion we obtain that $\frac{\pi}{4a}\frac{\varsigma'(1+\delta)}{\varsigma(1+\delta)}$ tends to $\frac{\pi}{4a}\left(-\frac{1}{\delta} + \gamma\right)$ while $\frac{\pi}{2}\frac{1}{a^2-(a-\delta)^2}$ tends to $\frac{\pi}{4a}\left(\frac{1}{\delta} + \frac{1}{2a}\right)$. This enables us to write the following equality: for $-2 < b < 1$,

$$\int_0^\infty \frac{t\arg(\varsigma(b+it))}{((1-b)^2+t^2)^2}dt = \frac{\pi}{4(1-b)}\left(\gamma - \frac{3}{2(1-b)}\right) +$$

$$\pi \sum_{\rho,\sigma_k>b,t_k>0}\left(\frac{1}{(1-b)^2+t_k^2} - \frac{(1-b)^2-(\sigma_k-b)^2+t_k^2}{((1-b)^2-(\sigma_k-b)^2+t_k^2)^2+4t_k^2(\sigma_k-b)^2}\right) \qquad (9)$$

and the following
THEOREM 2a. An equality

$$\int_0^\infty \frac{t\arg(\varsigma(b+it))}{((1-b)^2+t^2)^2}dt = \frac{\pi}{4(1-b)}\left(\gamma - \frac{3}{2(1-b)}\right) \qquad (10),$$

where $b$ is real $1/2 \leq b < 1$, holds true for some $b$ if and only if there are no Riemann function zeroes with $\sigma > b$. For $b=1/2$ this equality is equivalent to the Riemann hypothesis.

It is instructive to obtain the results of this Section in a different way. Let us first integrate $2\int_0^\infty \frac{t\ln(\varsigma(b+it))}{(a^2+t^2)^2}dt$ by parts:

$$2\int_0^\infty \frac{t\ln(\varsigma(b+it))}{(a^2+t^2)^2}dt = -\frac{1}{a^2+t^2}\ln(\varsigma(b+it))|_0^\infty + i\int_0^\infty \frac{1}{a^2+t^2}\frac{\varsigma'(b+it)}{\varsigma(b+it)}dt.$$ Thus we obtain

for the real part $\int_0^\infty \frac{2t\ln(|\varsigma(b+it)|)}{(a^2+t^2)^2}dt = \frac{\ln|\varsigma(b)|}{a^2} - \int_0^\infty \frac{1}{a^2+t^2}\text{Im}(\frac{\varsigma'(b+it)}{\varsigma(b+it)})dt$ and for imaginary part

$$\int_0^\infty \frac{2t\arg(\varsigma(b+it))}{(a^2+t^2)^2}dt = -\frac{\pi}{a^2} + \int_0^\infty \frac{1}{a^2+t^2}\text{Re}(\frac{\varsigma'(b+it)}{\varsigma(b+it)})dt \qquad (11)$$

for $b<1$ and $\int_0^\infty \frac{2t\arg(\varsigma(b+it))}{(a^2+t^2)^2}dt = \int_0^\infty \frac{1}{a^2+t^2}\text{Re}(\frac{\varsigma'(b+it)}{\varsigma(b+it)})dt$ for $b>1$. Let for a moment $1/2 < b < 1$, $a+b \neq 1$ and suppose that RH holds true. Then we can use the equality from Theorem 2 of our previous Note [1]
$\frac{a}{\pi}\int_{-\infty}^{\infty}\frac{\ln(|\varsigma(b+it)|)}{a^2+t^2}dt = \ln\left(\left|\frac{\varsigma(a+b)(a+b-1)}{a-b+1}\right|\right)$ and obtain by differentiating by $b$:

$$\frac{a}{\pi}\int_{-\infty}^{\infty}\frac{1}{a^2+t^2}\operatorname{Re}\frac{\varsigma'(b+it)}{\varsigma(b+it)}dt = \frac{\varsigma'(a+b)}{\varsigma(a+b)} + \frac{1}{a+b-1} - \frac{1}{a-b+1} \quad (12)$$

Substituting (12) into (11) we get the following equality

$$\frac{2a}{\pi}\int_0^{\infty}\frac{2t\arg(\varsigma(b+it))}{(a^2+t^2)^2}dt = -\frac{2}{a} + \frac{\varsigma'(a+b)}{\varsigma(a+b)} + \frac{1}{a+b-1} - \frac{1}{a-b+1} \quad (13).$$

In the general case, if $b<1$ we should add also the contributions of the Riemann function zeroes which now are given by the derivatives of $\frac{d}{db}\sum_{\sigma_k>b}\ln\left|\left(\frac{a+\sigma_k-b+it_k}{a-\sigma_k+b-it_k}\right)\right|$ by $b$; see [1]. Performing calculations and pairing the complex conjugate zeroes we get, of course, the equality (7).

REMARK 1. In view of $\int\frac{tdt}{(c^2+t^2)(d^2+t^2)} = \frac{1}{2(c^2-d^2)}\ln\frac{d^2+t^2}{c^2+t^2}$, similar approach apparently might be used to handle the integrals $\int_0^{\infty}\frac{t\arg(\varsigma(b+it))}{(c^2+t^2)(d^2+t^2)}dt$. However, as it follows from above, to do this we need first to analyze the integrals $\int_0^{\infty}\ln|\varsigma(b+it)|\ln\left(\frac{d^2+t^2}{c^2+t^2}\right)dt$, where it is not clear how to proceed. Hence for such a case only the direct approach taken above is appropriate.

2. **Integral equalities between integrals over the real axis and over the line** $z=b+it$

In recent paper [1] we showed that all equalities and theorems of the type considered can be reformulated to explore the function $\ln(\varsigma(z)\cdot(z-1))$ rather than $\ln(\varsigma(z))$. That is, for the case in hand we should consider integrals $\int_0^{\infty}\frac{t\arg(\varsigma(b+it)\cdot(b+it-1))}{(c^2+t^2)(d^2+t^2)}dt$, $\int_0^{\infty}\frac{t\arg(\varsigma(b+it)\cdot(b+it-1))}{(a^2+t^2)^2}dt$ and like. Such a modification kills the simple pole of the Riemann function occurring at $z=1$, hence the resulting relations are simply in certain aspects. We will not pursue this line of research in full in the current paper and will consider only one example, viz. the contour integral $\int_C\frac{-i(z-1/2)\ln(\varsigma(z)\cdot(z-1))}{z^2(z-1)^2}dz$ where in the contour $b=1/2$. On the left border of the contour the integral takes the form $-i\int_{-\infty}^{\infty}\frac{t\ln(\varsigma(1/2+it)\cdot(-1/2+it))}{(t^2+1/4)^2}dt$ hence it is a full analogue of integrals from eqs. (9, 10) when $b=1/2$. An analysis using the Laurent expansion of the Riemann function in the vicinity of the point $z=1$, $\ln(\varsigma(1+\delta)\cdot\delta)=\gamma\delta+...$, readily shows that in the interior of the contour now we have only one *simple* pole at $z=1$ with a residue equal to $-i\gamma/2$, hence Cauchy theorem gives $-i\int_C\frac{(z-1/2)\ln(\varsigma(z)\cdot(z-1))}{z^2(z-1)^2}dz=\pi\gamma$. Now, repeating [1] we introduce a "semi-contour" $C'''$ (see Fig. 1) which "semi-indents" the point $z=1$ thus obtaining for the imaginary part of the integrals:

$$-i\int_0^{\infty}\frac{t\ln(|\varsigma(1/2+it)\cdot(-1/2+it)|)}{(t^2+1/4)^2}dt-i\int_{1/2}^{\infty}\frac{(x-1/2)\ln(|\varsigma(x)\cdot(x-1)|)}{x^2(x-1)^2}dx=0,$$ that is

$$-\int_0^{\infty}\frac{t\ln(|\varsigma(1/2+it)\cdot(-1/2+it)|)}{(t^2+1/4)^2}dt=\int_{1/2}^{\infty}\frac{(x-1/2)\ln(|\varsigma(x)\cdot(x-1)|)}{x^2(x-1)^2}dx \quad (14).$$

Here the integral $\int_{1/2}^{\infty}\frac{(x-1/2)\ln(|\varsigma(x)\cdot(x-1)|)}{x^2(x-1)^2}dx$ is understood as a principal Cauchy value, that is

$$\int_{1/2}^{\infty}\frac{(x-1/2)\ln(|\varsigma(x)\cdot(x-1)|)}{x^2(x-1)^2}dx=\lim_{\varepsilon\to 0}(\int_{1/2}^{1-\varepsilon}\frac{(x-1/2)\ln(|\varsigma(x)\cdot(x-1)|)}{x^2(x-1)^2}dx+\int_{1+\varepsilon}^{\infty}\frac{(x-1/2)\ln(|\varsigma(x)\cdot(x-1)|)}{x^2(x-1)^2}dx).$$

Relation (14) contains integrals only of the logarithm of the modulus of the Riemann function, and thus is rather suitable for the numerical study, see below.

On the other hand, if RH is false, each Riemann zero $\sigma_k+it_k$ with $\sigma_k>1/2$ contributes

$$\text{Im}(-\pi n(\frac{1}{1/4-(\sigma_k-1/2+it_k)^2}-\frac{1}{1/4+t_k^2}))=\frac{-2\pi n(\sigma_k-1)t_k}{(1/4+t_k^2-(\sigma_k-1/2)^2)^2+4(\sigma_k-1/2)^2 t_k^2}$$

(see the analysis of the integral $\int_0^{\infty}\frac{t\arg(\varsigma(b+it))}{(a^2+t^2)^2}dt$ above) to the integral value. All these contributions for $t_k>0$ have the same sign, and hence (14) is equivalent to the Riemann hypothesis.

REMARK 2. This same hint with introducing the semi-contour *C'''*, first proposed by us in [1], can be repeated for other equalities of the paper thus leading to a number of relations similar to (14), that is connecting integrals of the logarithm of the modulus of the Riemann function along the real semi-axis and half of the vertical line *z=b+it*. We will not pursue this line of the research here.

## 4. Interpretation of above integrals as certain sums over the Riemann function zeroes. Analysis of Volchkov criterion

Our above equalities involve integrals containing the function $\arg(\varsigma(b+it))$ as a factor under the integral sign. In particular, for *b=1/2* we have a factor $\arg(\varsigma(1/2+it))$. But this is well known that this same factor is exactly the "non-trivial" part of an expression describing the number of the Riemann function zeroes lying in the critical strip and having the imaginary part *t<x*:

$$N(x) = 1 - \frac{x \ln \pi}{2\pi} + \frac{1}{\pi} \operatorname{Im} \ln \Gamma(\frac{1}{4} + \frac{ix}{2}) + \frac{1}{\pi} \arg(\varsigma(1/2 + ix)) \quad (15)$$

see e.g. p. 212 of Ref. [2] for details. This expression can be used for calculation of certain sums over the imaginary part of Riemann zeroes: the integral $\int_0^\infty G(x)dN(x)$, provided, of course, that it exists, is equal to the sum $\sum_{\rho, t_k > 0} G(t_k)$ taken over all Riemann zeroes. In practice, this calculation is interesting if the function *G(x)* is such that the integration by parts can be used: $\int_0^\infty G(x)dN(x) = G(x)N(x)\big|_0^\infty - \int_0^\infty g(x)N(x)dx$; here, evidently, $g(x) = dG(x)/dx$.

In 1995, Volchkov proposed to use this approach to check up the Riemann hypothesis [6]. He took $G(x) = \frac{1}{1/4 + x^2}$, hence $g(x) = -\frac{2x}{(1/4 + x^2)^2}$, and noted that the integral $\int_0^\infty G(x)dN(x) = -\int_0^\infty g(x)N(x)dx$, which by construction is equal to the sum $\Sigma = \sum_{\rho, t_k > 0} \frac{1}{1/4 + t_k^2}$ over Riemann zeroes, coincides with the real sum $\sum_\rho \frac{1}{\rho} = \sum_\rho \frac{1}{\sigma_k + it_k}$ over these zeroes if and only if the RH holds true. For the latter it is well known that (p. 82 Ref. [7]) $\sum_\rho \frac{1}{\rho} = \frac{\gamma}{2} + 1 - \frac{1}{2}\ln 4\pi$. Using this and calculating term by term the integral $\int_0^\infty g(x)N(x)dx$ applying (15), Volchkov then established that the equality

$$\frac{1}{\pi}\int_0^\infty \frac{2t\arg(\varsigma(1/2+it))}{(1/4+t^2)^2}dt = \gamma - 3 \qquad (16)$$

is equivalent to the RH. In this expression one immediately recognizes our equality (10) and Theorem 2a, of which Volchkov criterion is a particular case corresponding to $b=1/2$, see also [8].

REMARK 3. Volchkov's paper more or less explicitly contains a relation similar to (16) but actually he gave another form of his criterion. After obtaining (16), he considered the function $S_1(x) = \int_0^x \arg(\varsigma(1/2+it))dt$. In view of $S_1(x) = O(\ln x)$ (p. 214 of [2]), one can integrate by parts to obtain $\int_0^\infty \arg(\varsigma(1/2+ix))g(x)dx = -\int_0^\infty S_1(x)\frac{dg(x)}{dx}dx$. Then Volchkov introduced the relation $S_1(x) = \int_{1/2}^\infty (\ln|\varsigma(\sigma+ix)| - \ln|\varsigma(\sigma)|)d\sigma$ which follows from the Littlewood theorem concerning the integrals of logarithm of analytical functions taken round the rectangular contour with the vertices $(1/2, 1/2+iX, 1/2+X+iX, 1/2+X)$ for $X \to \infty$ and, noting that $\int_{1/2}^\infty \ln|\varsigma(\sigma)|d\sigma = Const$ which contributes nothing to the $\int_0^\infty S_1(x)g_1(x)dx$, obtained $\frac{32}{\pi}\int_0^\infty \frac{1-12t^2}{(1+4t^2)^3}dt\int_{1/2}^\infty \ln|\varsigma(\sigma+it)|d\sigma = 3-\gamma$. This is, of course, equivalent to (16).

We do not see much help in introducing the function $S_1(x) = \int_{1/2}^\infty (\ln|\varsigma(\sigma+ix)| - \ln|\varsigma(\sigma)|)d\sigma$ which apparently is difficult to handle, but would like to note that all our equalities and theorems can be modified to use integrals of this and similar functions instead of integrals of an argument of the Riemann zeta-function. As an example we give here the corresponding expression for our equality (6):

$$\int_0^\infty \frac{d}{dt}\left(\frac{t}{(9/4+t^2)(1/4+t^2)}\right)dt \int_{1/2}^\infty \ln|\varsigma(\sigma+it)|d\sigma = \frac{\pi}{4}\ln\left(\frac{27}{\pi^2}\right), \text{ that is}$$

$$\int_0^\infty \frac{t^4+t^2/2+9/16}{(9/4+t^2)^2(1/4+t^2)^2}dt \int_{1/2}^\infty \ln|\varsigma(\sigma+it)|d\sigma = \frac{\pi}{4}\ln\left(\frac{27}{\pi^2}\right) \qquad (17)$$

Again, (17) is equivalent to the Riemann hypothesis.

REMARK 4. Volchkov criterion can be generalized for other sums over the Riemann zeroes. For instance, one can take $\sum_\rho \frac{1}{s-\rho} = \sum_\rho \frac{1}{s-\sigma_k - it_k} = \sum_{\rho, t_k > 0} \frac{2(s-\sigma_k)}{(s-\sigma_k)^2 + t_k^2}$. For this sum we know that $\sum_\rho \frac{1}{s-\rho} = \frac{\varsigma'(s)}{\varsigma(s)} + \frac{1}{s-1} + \frac{1}{s} - \frac{1}{2}\ln\pi + \frac{1}{2}\psi\left(\frac{s}{2}\right)$ (pp. 80-82 [2]). On the other hand, on RH this sum is equal to the sum $\sum_{\rho, t_k > 0} \frac{2(s-1/2)}{(s-1/2)^2 + t_k^2}$ which can be calculated as an integral $\int_0^\infty G(x)dN(x) = -\int_0^\infty g(x)N(x)dx$ with $G(x) = \frac{2a}{a^2+x^2}$, $g(x) = G' = -\frac{4ax}{(a^2+x^2)^2}$; here $a = s-1/2$, $s \neq 1/2$, $N(x)$ is given by (15). One readily sees that the equality obtained in this way again will be equivalent to our Theorem 2 for $b=1/2$.

REMARK 5. In the end of the section 2.2 of the present paper we have shown that our equality (6) is a direct consequence of our Theorem 2 from the previous Note [1], i.e. of the equality $\frac{a}{\pi}\int_{-\infty}^\infty \frac{\ln(|\varsigma(b+it)|)}{a^2+t^2}dt = \ln\left|\frac{\varsigma(a+b)\cdot(a+b-1)}{a-b+1}\right| + \sum_{\sigma_k > b}\ln\left|\frac{a+\sigma_k-b+it_k}{a-\sigma_k+b-it_k}\right|$. It also has been demonstrated, respectively in the current and previous Notes, that for the particular case $b=1/2$ the equalities under discussion give correspondingly Volchkov's criterion $\frac{32}{\pi}\int_0^\infty \frac{1-12t^2}{(1+4t^2)^3}dt \int_{1/2}^\infty \ln|\varsigma(\sigma+it)|d\sigma = 3-\gamma$ [6] and Balazard-Saias-Yor's [9] criterion $\int_{-\infty}^\infty \frac{\ln(|\varsigma(1/2+it)|)}{1/4+t^2}dt = 0$ equivalent to the Riemann hypothesis. Thus it can be said that both these criteria are actually one and the same criterion differently written.

## 5. Numerical results.

To test some formulas given in this article we performed numerical calculations using the software packet *Mathematica*. More specifically, we used the function NIntegrate which permits as an input parameter the singularities of the integrand, in our case the zeroes of the zeta function.

First we investigated the criterion given by eq. (3). We integrated over the variable $t$ up to the value $t = 1000$, i.e. we considered the first 649 Riemann function zeroes [10]. The numerical convergence of the integral to $\frac{\pi}{20}\ln(18\pi^2/245)$ is quite fast, we obtained a difference of $3.21558907 \times 10^{-9}$. As a second example we took equation (6). Here the convergence to $\frac{\pi}{4}\ln(\pi^2/27)$ is of the same character, we obtained a difference of $3.21586984 \times 10^{-9}$ between the calculated and the theoretical values. A third calculation using equation (10) confirms that integrals of this kind behave numerically in the same manner. In fact, for $b=1/2$ the integral given by eq. (10) converges to $\frac{\pi}{2}(\gamma - 3)$ (of course provided the RH is true) and the difference is $3.21591683 \times 10^{-9}$.

Next we considered equation (10) (i.e. Theorem 2a) again but in another aspect. We could ask, how "good" this Theorem is, that is we may numerically investigate the dependence of the integral there on the parameter $\alpha$. To do these evaluations, we rewrite equation (10) as ($0 \leq \alpha < 1/2$):

$$\gamma(\alpha) = \frac{2-4\alpha}{\pi} \int_0^\infty \frac{t}{((1/2-\alpha)^2 + t^2)^2} \arg(\varsigma(1/2+\alpha+it))dt \qquad (18)$$

Now we can verify the convergence of the r.h.s. to the Euler constant γ as a function of α. The results are presented in Figure 2. They are in agreement with the theory: the approximation becomes better approaching the right critical strip's border $\alpha = 1/2$ where it is known that there are no zeroes.

The last numerical experiment involves equation (14). In this case we calculate the integrals until *t=5000* (i.e. for the l.h.s. we take into account the first 4520 zeroes of the Riemann function [10]). The r.h.s. doesn't pose any problem. We can calculate it with the method's parameter "PricipalValue" of the function "NIntegrate". We obtain for the l.h.s. of (10) the value 0.3946346584 (rounded at the $10^{th}$ digit) while the r.h.s. gives 0.3946344787, hence the difference is $3.607193 \times 10^{-7}$.

## 6. Conclusion

In this and our previous paper under the same title [1] we have established a number of criteria involving the integrals of the logarithm of the Riemann $\varsigma$-function and equivalent to the Riemann hypothesis. Our results include all earlier known criteria of this kind [3, 6, 9] which are certain particular cases of the general approach proposed.

Numerous other criteria of the type "$\int_{b-i\infty}^{b+i\infty} g(z)\ln(\varsigma(z))dz = f(b)$ is equivalent to RH" can be constructed following the lines of the present paper and [1], viz. selecting an appropriate function $g(z)$ and calculating the value of a contour integral exploiting the generalized Littlewood theorem. We would like to hope that among them *the* criterion which can be practically used to prove or disprove the RH might possibly appear.

REFERENCES

[1] S. K. Sekatskii, S. Beltraminelli, and D. Merlini, A few equalities involving integrals of the logarithm of the Riemann $\varsigma$-function and equivalent to the Riemann hypothesis, arXiv:0806.1596v1 [math.NT].

[2] E. C. Titchmarsh and E. R. Heath-Brown, The theory of the Riemann Zeta-function, Oxford, Clarendon Press, 1988.

[3] F. T. Wang, A note on the Riemann Zeta-function, Bull. Amer. Math. Soc. 52 (1946), 319.

[4] E. C. Titchmarsh, The theory of functions, Oxford, Oxford Univ. Press, 1939.

[5] J. van der Lune, H. J. J. Te Riele, and D. T. Winter, On the zeros of the Riemann Zeta function in the critical strip. IV., Math. Comp. 46, (1986) 667.

[6] V. V. Volchkov, On an equality equivalent to the Riemann hypothesis, Ukrainian Math. J., 47, (1995), 422.

[7] H. Davenport, Multiplicative number theory, Springer, New York, 2000.

[8] D. Merlini, The Riemann magneton of the primes, Chaos and complexity Lett., Nova Sci. Publ., New York, 2, (2006), 93.

[9] Balazard M., Saias E. and Yor M., Notes sur la fonction de Riemann, Advances in Mathematics, 143, (1999), 284.

[10] A. Odlyzko, http://dtc.umn.edu/~odlyzko/zeta_tables/zeros1

S. K. Sekatskii, Laboratoire de Physique de la Matière Vivante, IPSB, BSP 408, Ecole Polytechnique Fédérale de Lausanne, CH1015 Lausanne-Dorigny, Switzerland.

E-mail : serguei.sekatski@epfl.ch

S. Beltraminelli, CERFIM, Research Center for Mathematics and Physics, PO Box 1132, 6600 Locarno, Switzerland.

E-mail: Stefano.beltraminelli@ti.ch

D. Merlini, CERFIM, Research Center for Mathematics and Physics, PO Box 1132, 6600 Locarno, Switzerland.

E-mail: merlini@cerfim.ch

**Figures**

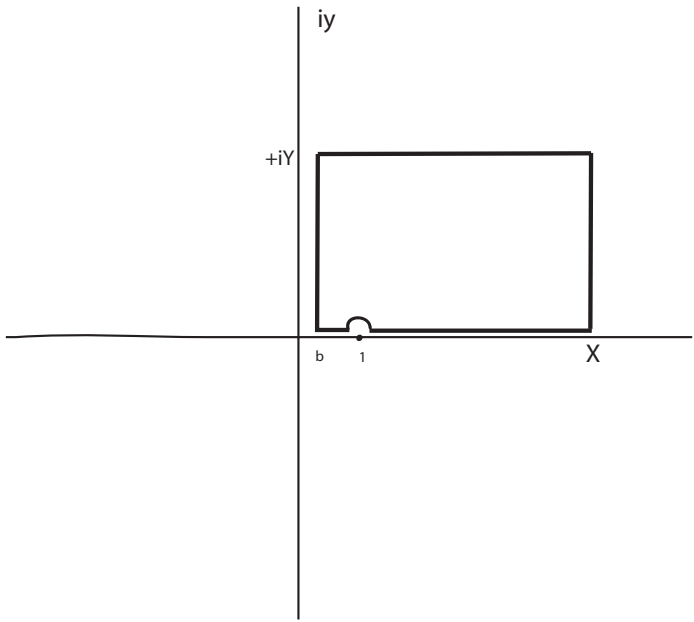

Figure 1. "Semi-contour" $C'''$ using for investigation of the integral in Section 3.

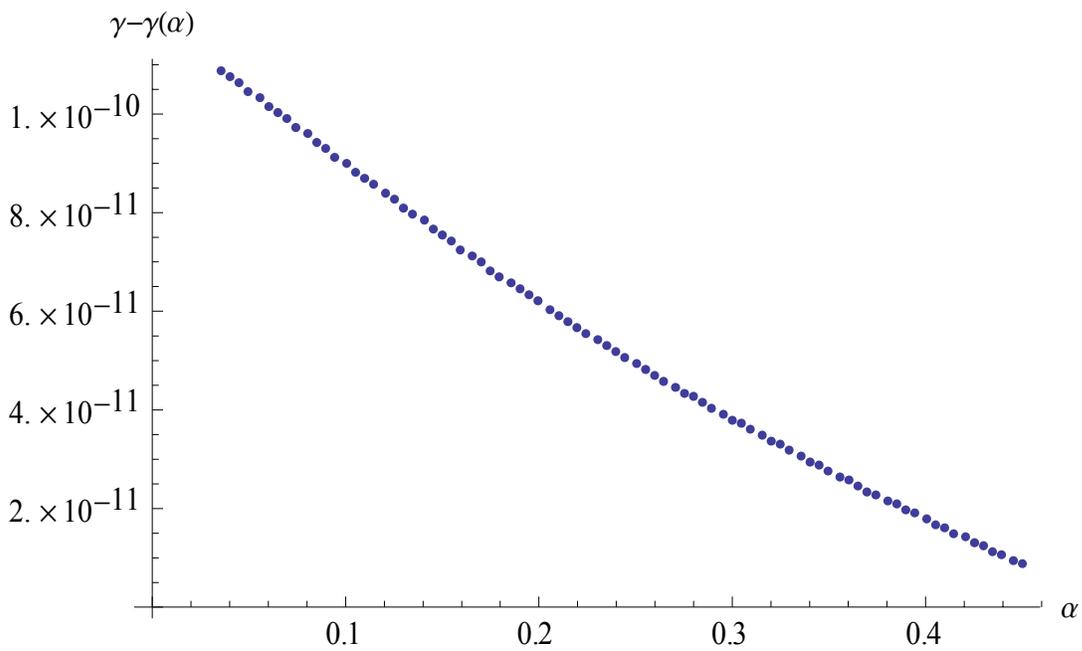

Figure 2. The convergence of the integral in eq. (19) to Euler constant $\gamma$ as a function of $\alpha$.